\documentclass[11pt]{amsart}
\usepackage[utf8]{inputenc}
\usepackage{amssymb, amsmath, amsthm, color, enumerate, mathabx, mathrsfs,latexsym}
\usepackage{bbm}

\usepackage{footmisc} 

\usepackage{pdfpages}

\usepackage{esint} 

\usepackage[normalem]{ulem}

\usepackage{tcolorbox}

\numberwithin{equation}{section}

\newcommand{\C}{\mathbb{C}}
\newcommand{\R}{\mathbb{R}}

\newcommand{\N}{\mathbb{N}}

\newcommand{\cD}{\mathcal{D}}

\newcommand{\cT}{\mathcal{T}}

\usepackage{accents}

\newcommand{\rank}{\mathrm{rank}}

\theoremstyle{plain}

\newtheorem{theorem}{Theorem}[section]

\newtheorem{proposition}[theorem]{Proposition}

\theoremstyle{definition}
\newtheorem{definition}[theorem]{Definition}

\newtheorem{remark}[theorem]{Remark}

\renewcommand{\Gamma}{\varGamma}
\renewcommand{\epsilon}{\varepsilon}

\renewcommand{\hat}{\widehat}
\renewcommand{\leq}{\leqslant}
\renewcommand{\geq}{\geqslant}

\author{Jennifer Duncan and Ana Vargas}
\title[codimension-2 restriction in $\R^5$]{Sharp $L^4$ restriction estimates for surfaces of codimension two in $\R^5$}
\begin{document}

\thanks{The authors were partially supported by grant  PID2022-142202NB-I00 AEI/10.13039/501100011033, J. Duncan is supported by the Juan de la Cierva grant JDC2023-050459-I, and the grant PID2024-158664NB-C2 financed by MICIU/ AEI / 10.13039/501100011033
/ FEDER, UE}
\begin{abstract}
We prove $L^p\rightarrow L^4$ Fourier adjoint restriction estimates for three-dimensional quadratic surfaces in $\R^5.$ These estimates are sharp up to the endpoint, and improve results by Oberlin \cite{O}, and Guo--Oh \cite{GO}.
\end{abstract}
\maketitle


\section{Introduction}

In this note we consider the three-dimensional submanifold
\begin{equation}\label{surface}
S=\{\Sigma(\xi):=(\xi,Q_A(\xi),Q_B(\xi))\,:\,\xi\in K\}\subset\R^5,
\end{equation}
where $Q_A(x):=x^\top A x$ and $Q_B(x):=x^\top B x$ are quadratic forms associated to symmetric matrices $A,B\in \R^{3\times 3}$, and $K\subset \R^3$ is a compact set. We define the associated extension operator $E_S$ to $S$ as, given $f\in L^1(K)$,
\begin{equation}\label{operator}
E_Sf(x):=\int_{K}e^{ix\cdot\Sigma(\xi)}f(\xi)d\xi.
\end{equation}

\

We will say that $R_S^*(p\rightarrow q)$ holds if there is a $C_{p,q}>0$ such that, for all $f\in L^p(K),$
\begin{equation}\label{estimate}
\|E_Sf\|_{L^q(\R^{5})}\le C_{p,q}\|f\|_{L^p(K)}.
\end{equation}

We say that a pair $(A,B)$ is irreducible if the matrices $A$ and $B$ are linearly independent and $\ker(A)\cap\ker(B)=\{0\}$, or equivalently, that there does not exist a linear change of coordinates in $\xi$ after which both $Q_A(\xi)$ and $Q_B(\xi)$ omit a common variable.

\

Two irreducible pairs $(A,B)$ and $(\tilde A,\tilde B)$ are \it equivalent\rm, and we denote this by $(A,B)\equiv (\tilde A, \tilde B)$, if there exist invertible matrices $M_1\in\R^{3\times 3}$ and $M_2\in\R^{2\times 2}$ such that
\begin{align}
    (\tilde A,\tilde B)^\top&=M_2\cdot(M^\top A M,M^\top BM)^\top\label{equiv}
\end{align}
where the action of $M_2$ on $(\R^{3\times 3})^2$ is the natural vector-valued extension of its action on $\R^2$, or alternatively, that
\begin{equation}
(Q_{\tilde A}(\xi), Q_{\tilde B}(\xi))^T=M_2\cdot(Q_A(M_1\xi), Q_B(M_1\xi))^T,
\end{equation}
for every $\xi\in\R^3.$ This means that, if we define
$$
\tilde S=\{(\xi, Q_{\tilde{A}}(\xi), Q_{\tilde{B}}(\xi))\,:\,\xi\in M_1^{-1}(K)\}\subset\R^5,
$$
then, after applying a change of co-ordinates, one may see that $R^*_S(p\rightarrow q)$ holds if and only if  $R^*_{\tilde S}(p\rightarrow q)$ holds. As a matter of convenience, we will say that a pair of quadratic forms $(Q_A,Q_B)$ is irreducible if $(A,B)$ is irreducible, and that $(Q_A,Q_B)$ is equivalent to $(Q_{\tilde{A}},Q_{\tilde{B}})$ if $(A,B)$ is equivalent to $(\tilde{A},\tilde{B})$. Note also that the irreducibility of the pair $(A,B)$ is preserved under transformations of the form \eqref{equiv}.

\

In \cite{GO} there is a classification of irreducible pairs $(A,B)$ that we shall now recall here. First of all, we have a condition introduced by M. Christ (\cite{Ch82}, \cite{Ch85}) and G. Mockenhaupt \cite{M}: a pair $(A,B)$ satisfies the (CM) condition if
 \begin{equation}\label{CM}
\int_{S^1}|\det(x_1A+x_2B)|^{-\gamma}d\sigma(x_1,x_2)\qquad\text{ for all }0<\gamma<\frac23.
\end{equation}
They showed that $(A,B)$ satisfies the (CM) condition if and only if $R^*(2\rightarrow 14/3)$ holds.
Christ also showed that for pairs satisfying the (CM) condition, $R^*(p\rightarrow q)$ fails for $q\le 10/3.$
D. Oberlin \cite{O} proved that for certain surfaces satisfying the (CM) condition $R^*(p\rightarrow q)$ holds for all $p>12/5$ and $q>4.$ He also found a necessary condition for those surfaces,
\begin{equation}\label{nec}
\frac3p+\frac7q\le 3
\end{equation}
The counterexample was given by $f=\chi_{[0,\epsilon]^3},$ for some $\epsilon>0,$ in which case
$|E_Sf|\ge\epsilon^3\chi_{[0,\epsilon^{-1}]^3\times[0,\epsilon^{-2}]^2}.$\par
Secondly, a pair $(A,B)$ satisfies the (D) condition if, for every $w\in \R^3\setminus\{0\}$, there exists a $\xi\in\R^3$ such that
\begin{equation}\label{D}
\det[A\xi,B\xi,w]\neq 0.
\end{equation}
And finally, for
$0\le d' \le 2,$ we say that the pair $(A,B)$ satisfies the ($R_{d'}$) condition if
\begin{equation}\label{Rd}
\min\limits_{hyperplane \;H\subset\R^3}
\max\limits_{\lambda_1,\lambda_2\in\R}
\rank( E_H^T(\lambda_1A+ \lambda_2B)E_H)=d',
\end{equation}
where $E_H\in\R^{3\times\dim(H)}$ is a matrix whose columns form a basis for $H$ (the exact choice of basis does not affect the value of the left-hand-side), The (D) condition and the ($R_d'$) condition were introduced in [DGS19] and
[GZ19b] to study decoupling inequalities for surfaces of codimension two.

\

In \cite{GO} it was proven that the (CM) condition is equivalent to the (D) together with the $(R_2$) condition. In \cite{GOR+} it was proven that condition (D) fails if and only if there is a non-trivial linear combination of $Q_A$ and $Q_B$ that is a complete square.

\

In \cite{GO} the author classified the surfaces that fail the (D) condition, and also those that satisfy the (D) condition plus the ($R_1$) or the ($R_0$) condition. Moreover they proved that for every irreducible pair ($P,Q$), $R^*_S(p\rightarrow q)$ holds for all $q>4$ and $p>\frac q{q-3}.$ They showed that, for surfaces failing the (D) condition and for surfaces satisfying the (D) and the ($R_0$) condition together, the result is sharp, up to the endpoints.

\

In this paper, we classify the surfaces satisfying the (CM) condition (i.e. the $(D)+(R_2)$) condition.

\begin{theorem}\label{classification}
If $(A,B)$ is an irreducible pair of matrices satisfying the (CM) condition, then, $(Q_A,Q_B)$ is equivalent to either $(\xi_2^2\pm  \xi_3^2, \xi_1^2\pm\xi_2^2 )$ or $( \xi_2\xi_3, \xi_1^2+\xi_2^2-\xi_3^2).$
\end{theorem}

 This result completes the classification of 3-dimensional quadratic surfaces in $\R^5$ initiated in \cite{GO}, which incidentally may be viewed as a simplified instance of the so-called `wild' problem of linear algebra, this being, given a field $\mathbb{K}$ and $n\in\N$, that of classifying the orbit spaces of pairs of matrices $(A,B)\in(\mathbb{K}^{n\times n})^{2}$ under the group action of $GL_n(\mathbb{K})$ through simultaneous conjugation, its namesake being due to its infamously intractible nature. For further reading on this topic we refer the reader to, for example, \cite{N, FKPS}. \par
 Moreover, we prove the following sharp restriction estimate.

\begin{theorem}\label{linearCM}
If $(A,B)$ satisfies the (CM) condition (i.e. $(A,B)$ satisfies the (D) and the ($R_2$) conditions), and $S$ is as defined in \eqref{surface}, then $R^*(p\rightarrow 4)$ holds for any $p>12/5.$
\end{theorem}
By interpolation with the trivial $R^*(p\rightarrow\infty)$ estimate, this theorem improves Oberlin's result, and furthermore, the necessary condition \eqref{nec} shows that Theorem \ref{linearCM} gives the best possible $L^4$-estimate, except for the endpoint. In \cite{GO}, the authors showed that all surfaces satisfying (D)+($R_1$)  conditions at the same time are equivalent to $(\xi_1\xi_2,\xi_1\xi_3+\xi_2^2),$ or to $(\xi_1\xi_2,\xi_1^2\pm\xi_3^2).$ For those surfaces, we have an estimate that improves upon theirs.

\begin{theorem}\label{linearDR1}
For those pairs $(A,B)$ satisfying both the $(D)$ and $(R_1)$  conditions, we instead have two different behaviours:
\begin{enumerate}
    \item If $(Q_A,Q_B)$ is equivalent to $(\xi_1\xi_2,\xi_1\xi_3+\xi_2^2),$ then $R^*_S(p\rightarrow 4)$ holds for any $p>3.$\\
    \item If $(Q_A,Q_B)$ is equivalent to $(\xi_1\xi_2,\xi_1^2\pm\xi_3^2),$ then $R^*_S(p\rightarrow 4)$ holds for any $p>8/3.$
\end{enumerate}
This result is sharp except, maybe, for the endpoint.
\end{theorem}
\begin{remark}By interpolation with the trivial $R^*(p\rightarrow\infty)$ estimate, this theorem improves Guo-Oh's result. Indeed, to see that this theorem gives the best possible estimates in $L^4$ except for the endpoint for $(D)+(R_1)$ surfaces.  we have the following examples:

\begin{enumerate}
    \item In the first surface, take $f=\chi_{[0,\epsilon]\times[0,\epsilon^{1/2}]\times[0,1]},$ for some $\epsilon>0.$ Then
$|E_Sf|\ge\epsilon^{3/2}\chi_{[0,\epsilon^{-1}]\times[0,\epsilon^{-1/2}]\times[0,1]
\times[0,\epsilon^{-3/2}]\times[0,\epsilon^{-1}]}.$ This gives the necessary condition
$$\frac4q+\frac3{2p}\le \frac32.$$
\item In the second surface, take $f=\chi_{[0,\epsilon]\times[0,1]\times[0,\epsilon]},$ for some $\epsilon>0.$ Then
$|E_Sf|\ge\epsilon^2\chi_{[0,\epsilon^{-1}]\times[0,1]\times[0,\epsilon^{-1}]
\times[0,\epsilon^{-1}]\times\times[0,\epsilon^{-2}]}.$ This gives the necessary condition
$$\frac5q+\frac2p\le 2.$$
\end{enumerate}

\end{remark}

\

Restriction estimates for submanifolds of codimension greater than one have been studied by many authors. Apart from those previously mentioned, A. D. Banner \cite{B} studied restriction estimates for certain surfaces not satisfying the (CM) condition; J. G. Bak and S. Lee \cite{BL} and Oberlin \cite{O02} obtained sharp Fourier restriction estimates for certain quadratic moment surfaces of high codimensions, and obtained sharp Fourier restriction estimates for certain quadratic moment surfaces of high codimensions.
J. K. Bak, J. Lee and S. Lee \cite{BLL} studied bilinear Fourier restriction estimates for surfaces of codimension larger than one and deduced a linear restriction estimate for the complex paraboloid in $\C^3.$ J. Lee and S. Lee \cite{LL} studied restriction estimates for complex hypersurfaces, which can also be viewed as special cases of real surfaces of codimension two. Notably, a restriction-type curvature condition for submanifolds of general codimension was introduced in work of Gressman \cite{G}, later utilised by Dendrinos, Musta\c t\v a, and Vitturi in their study of $2$-plane transform estimates related to fourier restriction \cite{DMV}.

\

In section \ref{classif} we prove Theorem \ref{classification}, in section \ref{CMt} we prove Theorem \ref{linearCM}, and in section \ref{DR1} we prove Theorem \ref{linearDR1}. Throughout, the notation $A\lesssim B$ signifies that there exists a constant $C>0$ such that $A\leq CB$. This constant is taken to be universal, in that it does not depend on any extraneous quantities or variables such as dimensions, functions, or exponents.

\section{Classification of (CM) surfaces}\label{classif}

The arguments in this section are very elementary and probably can be simplified. In the proof of Theorem \ref{classification}, we will make use, repeteadly of some basic observations:

\begin{remark}\label{basic}
a) Note that for any $c\in\R,$ $(A,B)\equiv(A,B-cA).$
To justify this assertion, just take
\begin{equation*}
M_2= \begin{pmatrix}
         1& 0&\\
         -c& 1
    \end{pmatrix}.
\end{equation*}

b) For any $\lambda\in \R,$ $(A,B)\equiv(A,\lambda B).$
To justify this, take $M_3=I,$ and
\begin{equation*}
M_2= \begin{pmatrix}
         1& 0\\
         0 & \lambda
    \end{pmatrix}.
\end{equation*}
\end{remark}

\

\noindent\it Proof of Theorem \ref{classification}\rm \quad We will begin by establishing two claims: \medskip\\
\noindent{\sl Claim 1 -- If $(A,B)$ is an irreducible pair that satisfies the (CM) condition, then either $\rank(A)>1$ or $\rank(A)>1.$}

Assume for contradiction that $\rank(A)\le 1,$ $\rank(B)\le 1.$, then $\ker(A)\cap\ker(B)\neq\{0\}$, and therefore the pair $(A,B)$ is reducible. In particular, it cannot satisfy \eqref{CM}.\medskip
\

\noindent{\sl Claim 2 -- Every irreducible pair $(A,B)$ is equivalent to another irreducible pair where the rank of one of the matrices is two.}

To prove this claim, note first that by Claim 1, we know that the range of one of the two matrices must be 2 or 3. Thus, to prove this claim, we only need to consider the case in which that rank is 3.

If $\rank A=3$ and $\rank B=1,$ take $\lambda_0,$ such that $\det(A-\lambda_0 B)=0$  (note that $\lambda_0$ is a real root of a cubic equation). Then, $(A,B)\equiv(A-\lambda_0 B,B)$ and, by Claim 1, $\rank(A-\lambda_0 B)=2.$

In case that $\rank(A)=\rank (B)=3,$ we can choose $\lambda_0\in\R$ such that $\det(A-\lambda_0 B)=0.$ Then,  $(A,B)\equiv(A-\lambda_0B, B)$ and $\rank(A-\lambda_0B)\le 2.$ If $\rank(A-\lambda_0B)= 2,$ then it satisfies the statement of the claim, while if $\rank(A-\lambda_0B)= 1,$ then, we are in the situation described in the previous paragraph and finish the proof as above.

\

\noindent{\sl Step 3}. We take $\rank(A)=2.$ Then, after a change of co-ordinates,
\begin{equation*}
A= \begin{pmatrix}
         0& 0&0 \\
         0& \alpha_2 & 0 \\
         0 & 0 & \alpha_3
    \end{pmatrix},
\end{equation*}
for some $\alpha_2\ne0,$ $\alpha_3\ne0.$
Hence, taking
\begin{equation*}
M_1= \begin{pmatrix}
         0& 0&0 \\
         0& \big(\sqrt{|\alpha_2|}\big)^{-1/2} & 0 \\
         0 & 0 & \big(\sqrt{|\alpha_3|}\big)^{-1/2}
    \end{pmatrix},
\end{equation*}
and
\begin{equation*}
M_2= \begin{pmatrix}
         \pm 1& 0\\
         0& 1
    \end{pmatrix},
\end{equation*}
we show that $(A,B)$ is equivalent to a pair where
\begin{equation}
A= \begin{pmatrix}
         0& 0&0 \\
         0& 1 & 0 \\
         0 & 0 & \pm 1
    \end{pmatrix}.
\end{equation}
We then consider the cases of different signs separately.\\
\noindent{\it Case 1.} Assume that
\begin{equation*}
A= \begin{pmatrix}
         0& 0&0 \\
         0& 1 & 0 \\
         0 & 0 & -1
    \end{pmatrix}.
\end{equation*}

Making an orthonormal change of variables, we change our pair for another equivalent pair where
\begin{equation}\label{canon}
A= \begin{pmatrix}
         0& 0&0 \\
         0& 0 & 1 \\
         0 & 1 & 0
    \end{pmatrix},
\end{equation}
and
\begin{equation*}
B= \begin{pmatrix}
         a& b&c \\
         b& d & f \\
         c & f &g
    \end{pmatrix}.
\end{equation*}

\noindent{\sl Claim.} $a\ne0.$
The reason is that, if $a=0,$ then,
$$
\det(x_1A+x_2B)=x_2^2R(x_1,x_2),
$$
for some polynomial $R$, and thus by checking the definition directly, one then quickly sees that $(A,B)$ cannot satisfy the (CM) condition.

Since $a\ne0,$ by Remark \ref{basic}, we can assume that $a>0.$ Then,
$$Q_A(\xi)=2\xi_2\xi_3,$$
and
\begin{align*}
    Q_B(\xi)&=a\xi_1^2+d\xi_2^2+g\xi_3^2+2b\xi_1\xi_2+2 c\xi_1\xi_3+2f\xi_2\xi_3\\
    &=a\big(\xi_1+b\xi_2/a+c\xi_3/a)^2+d'\xi_2^2+g'\xi_3^2+
2f'\xi_2\xi_3,
\end{align*}
where $d':=d-b^2a^{-1}$, $g':=g-c^2a^{-1}$ and $f':=f-bca^{-1}$. If we now take
\begin{equation*}
M_1= \begin{pmatrix}
         1& -b/a& -c/a \\
         0& 1 & 0 \\
         0 & 0 &1
    \end{pmatrix},
\end{equation*}
we replace our pair by another equivalent where $P$ is as in \eqref{canon} and
\begin{equation*}
B= \begin{pmatrix}
         a& 0& 0 \\
         0& d' & f' \\
         0 & f' &g'
    \end{pmatrix}.
\end{equation*}
By Remark \ref{basic}, we can then replace $B$ with
\begin{equation*}
\begin{pmatrix}
         a& 0& 0 \\
         0& d' & 0 \\
         0 & 0 &g'
    \end{pmatrix}.
\end{equation*}
If $d'=0$ or $g'=0,$ then
$$
\det(x_1A+x_2B)=ax_2x_1^2,
$$
which as we have remarked earlier cannot satisfy the (CM) condition. Hence $d'\ne0$ and $g'\ne0.$
Taking
\begin{equation*}
M_1= \begin{pmatrix}
         1/\sqrt{a}& 0& 0 \\
         0& 1/\sqrt{|d'|} & 0 \\
         0 & 0 &1/\sqrt{|g'|}
    \end{pmatrix},
\end{equation*}
and
\begin{equation*}
M_2= \begin{pmatrix}
         \alpha& 0 \\
         0& 1
    \end{pmatrix},
\end{equation*}
for an appropriate $\alpha,$ we arrive at a pair where
\begin{equation*}
B= \begin{pmatrix}
         1& 0& 0 \\
         0&\pm1 & 0 \\
         0 & 0 &\pm 1
    \end{pmatrix},
\end{equation*}
and $A$ as in \eqref{canon}.

We analyze now the different possibilities.
If
\begin{equation*}
B= \begin{pmatrix}
         1& 0& 0 \\
         0&1 & 0 \\
         0 & 0 & -1
    \end{pmatrix},
\end{equation*}
we have finished.

\

If
\begin{equation*}
B= \begin{pmatrix}
         1& 0& 0 \\
         0&1 & 0 \\
         0 & 0 & 1
    \end{pmatrix},
\end{equation*}
then, taking
\begin{equation*}
M_1= \begin{pmatrix}
         1& 0& 0 \\
         0&\frac1{\sqrt2} & -\frac1{\sqrt2} \\
         0 & \frac1{\sqrt2} & \frac1{\sqrt2}
    \end{pmatrix},
\end{equation*}
we show that the pair $(A,B)$ is equivalent to $(\tilde{A},\tilde{B})$ where
\begin{equation*}
\tilde A= \begin{pmatrix}
         0& 0& 0 \\
         0& 1 & 0 \\
         0 & 0 &- 1
    \end{pmatrix},
\end{equation*}
\begin{equation*}
\tilde B= \begin{pmatrix}
         1& 0& 0 \\
         0&1 & 0 \\
         0 & 0 & 1
    \end{pmatrix}.
\end{equation*}
and therefore equivalent to $(\tilde A,\hat B),$ where
\begin{equation*}
\hat B:= \tilde B+\tilde A=\begin{pmatrix}
         1& 0& 0 \\
         0&2& 0 \\
         0 & 0 & 0
    \end{pmatrix},
\end{equation*}
which gives us the expected result, after application of Remark \ref{basic} with
\begin{equation*}
M_2= \begin{pmatrix}
         1& 0 \\
         0& 1/2
    \end{pmatrix}
\end{equation*}
and
\begin{equation*}
M_1= \begin{pmatrix}
          \sqrt2& 0& 0 \\
         0& 1 & 0 \\
         0 & 0 & 1
    \end{pmatrix}
\end{equation*}
Finally, if
\begin{equation*}
B= \begin{pmatrix}
         1& 0& 0 \\
         0&-1 & 0 \\
         0 & 0 & -1
    \end{pmatrix},
\end{equation*}
we proceed similarly.

\

\noindent{\it Case 2.} Assume now that
\begin{equation*}
A= \begin{pmatrix}
         0& 0&0 \\
         0& 1 & 0 \\
         0 & 0 & 1
    \end{pmatrix}
\end{equation*}
and
\begin{equation*}
B= \begin{pmatrix}
         a& b&c \\
         b& d & f \\
         c & f &g
    \end{pmatrix}.
\end{equation*}
As before, we can assume that $a>0.$ Then, we take
\begin{equation*}
M_1= \frac1{\sqrt{a}}\begin{pmatrix}
         1& -b/a&-c/a \\
         0& 1& 0 \\
         0 & 0 & 1
    \end{pmatrix}.
\end{equation*}
This shows that $(A,B)\equiv(A,\tilde B),$ where
\begin{equation*}
\tilde B= M_1^\top BM_1=\begin{pmatrix}
         1&0&0\\
         0& d'& f' \\
         0 & f' & g'
    \end{pmatrix}.
\end{equation*}
Then, there exists an orthogonal matrix
\begin{equation*}
\tilde M_1= \begin{pmatrix}
         1& 0&0 \\
         0& \nu& \rho \\
         0 & -\rho & \nu
    \end{pmatrix}
\end{equation*}
such that
\begin{equation*}
\hat B=\tilde M_1^\top\tilde B\tilde M_1= \begin{pmatrix}
         1&0&0 \\
         0& d''& 0 \\
         0 & 0 & g''
    \end{pmatrix},
\end{equation*}
for some $d'',g''\in\R$, and we note that $A=\tilde M_1^\top A\tilde M_1$. It follows from Remark \ref{basic} that $(A,B)\equiv(A,\hat B)\equiv(A,\hat B-g''A),$ where
\begin{equation*}
\hat B-g''A=\begin{pmatrix}
         1&0&0 \\
         0& d'''& 0 \\
         0 & 0 & 0
    \end{pmatrix}.
\end{equation*}
Note that by the (CM) condition, $d'''\ne0.$
Application of Claim \ref{basic} with
\begin{equation*}
\tilde M_1= \begin{pmatrix}
         \sqrt{|d''|}& 0&0 \\
         0& 1 &0 \\
         0 & 0 & 1
    \end{pmatrix}
\end{equation*}
and
\begin{equation*}
\tilde M_2= \begin{pmatrix}
         1& 0\\
         0& \frac1{|d''|}
    \end{pmatrix}
\end{equation*}
finishes the proof.

\qed


\section{$L^4$ estimate for the (CM) case}\label{CMt}
The arguments in this and next section are classical. To prove Theorem \ref{linearCM}, we will use a bilinear estimate, and as such we will consider pairs of rectangles in $\R^3$ satisfying a \it separation \rm condition.

\begin{definition}
Let $j,k,\ell\in\N$, and let $\mathcal D_{j,k,\ell}$ denote the family of dyadic rectangles of dimension $2^{-j}\times2^{-k}\times 2^{-\ell}$ in $[0,1]^3.$ We say that two rectangles $\tau,\tau'\in\mathcal D_{j,k,\ell}$ are \it separated\rm, and indicate that by $\tau\sim\tau',$ if
$|\xi_1-\eta_1|\sim 2^{-j}$, $|\xi_2-\eta_2|\sim 2^{-k}$, $|\xi_3-\eta_3|\sim 2^{-\ell}$ for all $\xi\in\tau$ and $\eta\in\tau'.$
\end{definition}
The following bilinear estimate is our key proposition, from which Theorem \ref{CM} quickly follows
\begin{proposition}\label{bilinearCM} Let $\tau, \tau'\in\mathcal D_{j,k,\ell}$ be separated rectangles in $[0,1]^3$, and let $f$ and $g$ be functions supported in $\tau$ and $\tau'$ respectively. Let  $(Q_A(\xi),Q_B(\xi))=(\xi_2^2\pm  \xi_3^2, \xi_1^2\pm\xi_2^2 ),$ or $(Q_A(\xi),Q_B(\xi))=( \xi_2\xi_3, \xi_1^2+\xi_2^2-\xi_3^2).$ Suppose that $S$ is the quadratic 3-surface determined by $(Q_A,Q_B)$ as in \eqref{surface}, then,
$$
\|E_S f E_Sg\|_{L^2}\le C 2^{\frac{j+k+\ell}6}\|f\|_{L^2}\|g\|_{L^2}.
$$
\end{proposition}

\

We postpone the proof of this proposition to the end of the section.

\

\noindent\it Proof of Theorem \ref{linearCM}.\rm

By Theorem \ref{classification} we can assume that $(Q_A,Q_B)$ is a pair as in the statement of Proposition \ref{bilinearCM}. By rescaling, we can assume that $K=[0,1]^3.$ The proof of Theorem \ref{linearCM} follows the ideas of the \it bilinear implies linear \rm argument of \cite{TVV}. Consider the case of $f=\chi_\Omega$ for some $\Omega\subset [0,1]^3$.
Letting $\Delta\subset[0,1]^3\times[0,1]^3$ denote the diagonal, we take a Whitney decomposition of $[0,1]^3\times[0,1]^3\setminus\Delta$, yielding a cover $\cT$ of $([0,1]^3\times[0,1]^3)\setminus\Delta$ consisting of pairs of caps $(\tau,\tau')$ such that, for some $j,k,\ell\in\N_0$, $\tau,\tau'\in\cD_{j,k,\ell}$ and $\tau\sim\tau'$. Let $\cT_{j,k,\ell}:=\cT\cap \cD_{j,k,\ell}^2$. Then,
\begin{align*}
    \|E_S\chi_\Omega\|_{L^4}^2&=\|E_S\chi_\Omega E_S\chi_\Omega\|_{L^2}\\
    &=\bigg\|\sum_{j,k,\ell\geq 0}\sum_{(\tau,\tau')\in\mathcal T_{j,k,\ell}} E_S\chi_{\Omega\cap\tau} E_S\chi_{\Omega\cap\tau'} \bigg\|_{L^2}\\
    &\leq\sum_{j,k,\ell\geq 0}\bigg\|\sum_{(\tau,\tau')\in\mathcal T_{j,k,\ell}} E_S\chi_{\Omega\cap\tau} E_S\chi_{\Omega\cap\tau'} \bigg\|_{L^2}\\
    &\lesssim \sum_{j,k,\ell\geq 0}\bigg(\sum_{(\tau,\tau')\in\mathcal T_{j,k,\ell}} \|E_S\chi_{\Omega\cap\tau} E_S\chi_{\Omega\cap\tau'}\|_{L^2}^2\bigg)^{1/2},
\end{align*}
where we have applied the Minkowski integral inequality and the second inequality is a consequence of the almost disjointness of the family of sets $\{\tau+\tau'\;|\;\tau,\tau'\in\mathcal D_{j,k,\ell},\;\; \tau\sim\tau'\}.$

\

We note that, by Proposition \ref{bilinearCM}, for fixed $j,k,\ell\in\N_0$
\begin{align*}
\sum_{(\tau,\tau')\in\cT_{j,k,\ell}} \|E_S\chi_{\Omega\cap\tau} E_S\chi_{\Omega\cap\tau'}\|_{L^2}^2&\le 2^{\frac{j+k+\ell}3}\sum_{(\tau,\tau')\in\mathcal T_{j,k,\ell}}|\Omega\cap\tau||\Omega\cap\tau'|\\
&\le 2^{\frac{j+k+\ell}3}\sum_{\tau\in\mathcal D_{j,k,\ell}}|\Omega\cap\tau|^2\\
&\le 2^{\frac{j+k+\ell}3}\sum_{\tau\in\mathcal D_{j,k,\ell}} \min\{|\Omega|, 2^{-(j+k+\ell)}\}|\Omega\cap\tau|\\
&=2^{\frac{j+k+\ell}3}\min\{|\Omega|, 2^{-(j+k+\ell)}\} |\Omega|.
\end{align*}
Hence,
\begin{align*}
  \|E_S\chi_\Omega\|_{L^4}^2&=\sum_{j,k,\ell\geq 0}\big(2^{\frac{j+k+\ell}3}\min\{|\Omega|, 2^{-(j+k+\ell)}\} |\Omega|\big)^{1/2}\\
    &\le  \sum_{\substack{j,k,\ell\geq 0\\ 2^{-(j+k+\ell)}\ge|\Omega|}}2^{(j+k+\ell)/6}|\Omega|+ \sum_{\substack{j,k,\ell\geq 0\\ 2^{-(j+k+\ell)}<|\Omega|}}2^{-(j+k+\ell)/3}|\Omega|^{1/2}.
\end{align*}

Let $\epsilon>0$ be such that $\frac1p=\frac5{12}-\epsilon.$ Then,
\begin{align*}
 &=  \sum_{j,k,\ell \ge 0}2^{-2\epsilon(j+k+\ell)}|\Omega|^{-\frac16+1-2\epsilon}+ \sum_{j,k,\ell\geq 0}2^{-2\epsilon(j+k+\ell)}|\Omega|^{\frac13+\frac12-2\epsilon}\\
 &=|\Omega|^{\frac56-2\epsilon}=|\Omega|^{2(\frac5{12}-\epsilon)}=|\Omega|^{\frac2p}=\|\chi_\Omega\|_{L^p}^2.
\end{align*}
\qed

\

\noindent\it Proof of Proposition \ref{bilinearCM}.\rm

Here we shall follow a classical $L^2$ approach based in a change of variables, similar to those used by L. Carleson and P. Sj\"olin \cite{CS}, C. Fefferman \cite{F} and A. Zygmund \cite{Z}.
Let $\tau,\tau'\in \cD_{j,k,\ell}$ be separated rectangles, and write $\tau'=I_1\times I_2\times I_3,$ where $|I_1|=2^{-j}$, $|I_2|=2^{-k}$, and $|I_3|=2^{-\ell}.$ We define a function $\Phi:[0,1]^6\rightarrow\R^5$ via the mapping
\begin{equation*}
    \Phi(\xi,\eta)
    :=(\xi_1-\eta_1,\xi_2-\eta_2,\xi_3-\eta_3,Q_A(\xi)-Q_A(\eta),Q_B(\xi)-Q_B(\eta)).
\end{equation*}
and, given a $1\leq j\leq 3$ and $\eta_j\in[0,1]$, we define $\Phi_{j,\eta_j}:[0,1]^5\rightarrow\R^5$, $\Phi_{j,\eta_j}(\xi,\widehat{\eta}):=\Phi(\xi,\eta)$, where $\widehat{\eta}\in[0,1]^2$ denotes $\eta=(\eta_1,\eta_2,\eta_3)$ with the $j^{th}$ variable omitted. We may view $\Phi_{j,\eta_j}$ as a `frozen' version of the function $\Phi$ where we have fixed the $j^{th}$ variable at the value $\eta_j\in[0,1]$.

Fixing $\eta_3\in I_3$, the Jacobian of $\Phi_{3,\eta_3}$ is given by
\begin{equation*}
    d\Phi_{3,\eta_3}(\xi,\eta_1,\eta_2):=\begin{pmatrix}
        1 & 0 & 0 & -1 & 0\\
        0 & 1 & 0 & 0 & -1\\
        0 & 0 & 1 & 0 & 0\\
        \partial_1Q_A(\xi) & \partial_2Q_A(\xi) & \partial_3Q_A(\xi) & -\partial_1Q_A(\eta) & -\partial_2Q_A(\eta)\\
        \partial_1Q_B(\xi) & \partial_2Q_B(\xi) & \partial_3Q_B(\xi) & -\partial_1Q_B(\eta) & -\partial_2Q_B(\eta)
    \end{pmatrix}.
\end{equation*}
Via elementary row-column operations,
\begin{align*}
    \det(d\Phi_{3,\eta_3}(\xi,\eta_1,\eta_2))=&\\
    &\hspace{-4cm}\begin{vmatrix}
         1 & 0 & 0 & 0 & 0\\
         0 & 1 & 0 & 0 & 0\\
         0 & 0 & 1 & 0 & 0\\
         \partial_1Q_A(\xi) &\partial_2Q_A(\xi) &\partial_3Q_A(\xi) &\partial_1Q_A(\xi)-\partial_1Q_A(\eta) &\partial_2Q_A(\xi)-\partial_2Q_A(\eta)\\
        \partial_1Q_B(\xi) & \partial_2Q_B(\xi) & \partial_3Q_B(\xi) & \partial_1Q_B(\xi)-\partial_1Q_B(\eta) & \partial_2Q_B(\xi)-\partial_2Q_B(\eta)
    \end{vmatrix}.
\end{align*}
Hence the determinant of this matrix is given by the $2\times 2$ determinant
\begin{equation}\label{determinant}
  \det(d\Phi_{3,\eta_3}(\xi,\eta_1,\eta_2))=\left|\begin{matrix}
         \partial_1Q_A(\xi)-\partial_1Q_A(\eta) &\partial_2Q_A(\xi)-\partial_2Q_A(\eta)\\
            \partial_1Q_B(\xi)-\partial_1Q_B(\eta) & \partial_2Q_B(\xi)-\partial_2Q_B(\eta)
    \end{matrix}\right|
\end{equation}
Now, after some elementary calculations, we see that $\Phi_{3,\eta_3}$ is a $\mathcal C^\infty$ injective immersion on $\tau\times I_1\times I_2$ in both the cases that $(Q_A,Q_B)=(\xi_1^2\pm\xi_2^2,\xi_2^2\pm\xi_3^2)$ and that $(Q_A,Q_B)=(\xi_1\xi_2,\xi_1^2+\xi_2^2-\xi_3^2). $ The inverse function theorem shows that it is a $\mathcal C^\infty$ diffeomorphism. We may thus apply a change of variables via $\Phi_{3,\eta_3}$ to the product $E_Sf(x)\overline{E_Sg}(x)$, after applying Fubini to isolate the variable $\eta_3$.
\begin{align*}
&E_Sf(x_1,x_2,x_3,x_4,x_5)\overline{E_Sg(x_1,x_2,x_3,x_4,x_5)}\\&=\iint_{\tau\times\tau'} f(\xi)\overline{g(\eta)}e^{2\pi i (x_1(\xi_1-\eta_1)+x_2(\xi_2-\eta_2)+x_3(\xi_3-\eta_3)+x_4(P(\xi)-P(\eta))+
x_5(Q(\xi_1)-Q(\eta))}\,d\xi d\eta\\
&=\int_{I_3}\bigg[\iint_{\Phi_{3,\eta_3}(\tau\times I_1\times I_2)} f(\Phi^{-1}_{3,\eta_3}(u)) \overline{g(\Phi^{-1}_{3,\eta_3}(u))}|d\Phi_{3,\eta_3}^{-1}(u)|e^{2\pi i x\cdot u}\,d u\bigg]d\eta_3\\
&=\int_{I_3}\widehat{F_{\eta_3}}(x)\,d\eta_3,
\end{align*}
where, letting $\chi_{\Phi_{3,\eta_3}(\tau\times I_1\times I_2)}$ denote the characteristic function of the image $\Phi_{3,\eta_3}(\tau\times I_1\times I_2)$,
$$
F_{\eta_3}(u):=f\otimes \widebar{g}(\Phi^{-1}_{3,\eta_3}(u))|d\Phi_{3,\eta_3}^{-1}(u)|\chi_{\Phi_{3,\eta_3}(\tau\times I_1\times I_2)}.
$$
Then, by Minkowski's inequality, the Cauchy-Schwarz inequality, and Plancherel's theorem,
\begin{align*}
    \|E_SfE_Sg&\|_{L^2}\le\int_{I_3}\|\widehat{F_{\eta_3}}\|_{L^2}\,d\eta_3\\
    &\leq 2^{-\ell/2}\bigg(\int_{I_3}\|\widehat{F_{\eta_3}}\|_{L^2}^2d\eta_3\bigg)^{1/2}\\
    &=2^{-\ell/2}\bigg(\int_{I_3}\|F_{\eta_3}\|_{L^2}^2d\eta_3\bigg)^{1/2}\\
    &=2^{-\ell/2}\bigg(\int_{I_3}\iint_{\Phi_{3,\eta_3}(\tau\times I_1\times I_2)} |f\otimes g(\Phi^{-1}_{3,\eta_3}(u))|^2|d\Phi_{3,\eta_3}^{-1}(u)|^{-2}
    \,du\,d\eta_3\bigg)^{1/2}.
\end{align*}
Changing back to our original variables, and noting that $d\Phi_{3,\eta_3}^{-1}(u)=d\Phi_{3,\eta_3}(\Phi_{3,\eta_3}^{-1}(u))$, we obtain the key estimate
\begin{align}\label{integral}
    &= 2^{-\ell/2}\bigg(\int_\tau\int_{\tau'}
    \frac{f(\xi)^2g(\eta)^2}{|d\Phi_{3,\eta_3}(\xi,\eta_1,\eta_2)|}d\xi d\eta\bigg)^{1/2}.
\end{align}

From this point onwards we consider the different cases separately.
\begin{itemize}
    \item In the case
$(Q_A(\xi),Q_B(\xi))=(\xi_1^2\pm\xi_2^2,  \xi_2^2\pm\xi_3^2),$ we compute $|d\Phi_{3,\eta_3}(\xi,\eta_1,\eta_2)|=|\xi_1-\eta_1||\xi_2-\eta_2|\sim 2^{-j-k}.$
Thus, by \eqref{determinant} and \eqref{integral},
\begin{align}
\|E_SfE_Sg\|_{L^2}\lesssim 2^{-\frac\ell2}2^{\frac{j+k}2}\|f\|_{L^2}\|g\|_{L^2}
\end{align}
Assume first that $\ell\ge j, k$, in which case $\frac{-\ell+j+k}2\le \frac{j+k+\ell}6,$ and so we are done. If $j\ge j, \ell$ or $k\ge j, \ell$, a similar argument, freezing $\eta_1$ or $\eta_2$ instead of $\eta_3$ gives the result.
\item In the case $(Q_A(\xi),Q_B(\xi))=(\xi_2\xi_3, \xi_1^2+\xi_2^2-\xi_3^2),$ then,
$$2^{-\ell/2}|d\Phi_{3,\eta_3}(\xi,\eta_1,\eta_2)|^{-1/2}\sim 2^{j/2}\le 2^{(j+k+\ell)/6}
$$
when $j\le k$ and $j\le \ell,$ while
\begin{align*}
2^{-j/2}|d\Phi_{1,\eta_1}(\xi,\eta_1,\eta_2)|^{-1/2}&\sim 2^{-j/2}\big(|\xi_2-\eta_2|^2+|\xi_3-\eta_3|^2\big)^{-1/2}\\
&\le
\min\{2^{-j/2}2^{k},2^{-j/2}2^{\ell}\}\le 2^{(j+k+\ell)/6},
\end{align*}
when $j\ge k$ or $j\ge \ell.$ So, by repeating the change of variables argument, this time freezing $\eta_1$ instead of $\eta_3$, we obtain the same estimate. One should bear in mind the small technical detail that here we have a failure of injectivity for $\Phi_{1,\eta_1}$ along the quadratic hypersurface
$$Z_1:=\{(\xi_1,\xi_2,\xi_3;\eta_1,\eta_2)\in\R^5:|\xi_2-\eta_2|^2=|\xi_3-\eta_3|^2\}$$
however since this is the zero-set of a non-zero polynomial, it is closed and null, so we may simply apply the change of variables to the set $(\tau\times I_2\times I_3)\setminus Z_1$ instead without affecting the rest of the argument.
\end{itemize}

\hfill\qed


\section{$L^4$ estimate for the (D)+$(R_1)$ case}\label{DR1}
We can prove Theorem \ref{linearDR1} by following an argument similar to that in the previous section, using the following bilinear estimates:

\begin{proposition}\label{bilinearDR1} Let $\tau, \tau'\in\mathcal D_{j,k,\ell}$ be separated rectangles in $[0,1]^3.$ Let $f$ and $g$ functions supported in $\tau$ and $\tau'$ respectively.

\begin{enumerate}
    \item If $(Q_A(\xi),Q_B(\xi))=(\xi_1\xi_2,\xi_1^2\pm\xi_3^2),$ then
$$
\|E_S f E_Sg\|_{L^2}\lesssim 2^{\frac{j+k+\ell}4}\|f\|_{L^2}\|g\|_{L^2}.
$$
\item If $(Q_A(\xi),Q_B(\xi))=(\xi_1\xi_2,\xi_1\xi_3+\xi_2^2),$ then
$$
\|E_S f E_Sg\|_{L^2}\lesssim 2^{\frac{j+k+\ell}3}\|f\|_{L^2}\|g\|_{L^2}.
$$
\end{enumerate}

\end{proposition}

\

\noindent\it Proof of Proposition \ref{bilinearDR1}.\rm

Here we repeat the argument of the proof of Proposition \ref{bilinearCM}. First of all, we note that, in the case that $(Q_A,Q_B)=(\xi_1\xi_2,\xi_1^2\pm\xi_3^2)$, we find that $\Phi_{1,\eta_1}$ and $\Phi_{3,\eta_3}$ act as diffeomorphisms on $\tau\times I_2\times I_3$ and $\tau\times I_1\times I_2$ respectively for all $\eta_1,\eta_3\in I_3$, however, in the case that $(Q_A,Q_B)=(\xi_1\xi_2,\xi_1\xi_3+\xi_2^2)$, while $\Phi_{1,\eta_1}$ is a diffeomorphism on $\tau\times I_2\times I_3$ for all $\eta_1\in I_1$, $\Phi_{3,\eta_3}$ on the otherhand is only injective on $\tau\times I_1\times I_2$ away from the quadratic hypersurface
$$Z_3:=\{(\xi_1,\xi_2,\xi_3;\eta_1,\eta_2)\in\R^5:(\xi_1-\eta_1)(\xi_3-\eta_3)=2|\xi_2-\eta_2|^2\}.$$
As before, this technicality does not meaningfully affect any of the calculations involved. Carrying out the argument again we obtain \eqref{integral} for this new collection of phases.\par\begin{enumerate}
    \item If $(Q_A(\xi),Q_B(\xi))=(\xi_1\xi_2,\xi_1^2\pm\xi_3^2),$ the Jacobian computations give us the following:
\begin{itemize}
    \item If $\ell\ge j,$ then freezing the third variable we obtain
\begin{align*}
    2^{-\ell/2}|d\Phi_{3,\eta_3}(\xi,\eta_1,\eta_2)|^{-1/2}&=
2^{-\ell/2}|\xi_1-\eta_1|^{-1}\\
&\sim2^{j-\ell/2}\le 2^{j/2}\le 2^{\frac{j+\ell}4}
\le 2^{\frac{j+k+\ell}4}.
\end{align*}
which together with \eqref{integral} gives us the desired estimate.
\item If $\ell\le j,$ then freezing the first variable this time
\begin{align*}
    2^{-j/2}|d\Phi_{1,\eta_1}(\xi,\eta_1,\eta_2)|^{-1/2}&=
2^{-j/2}\big(|\xi_1-\eta_1||\xi_3-\eta_3|\big)^{-1/2}\\
&\sim2^{-j/2+j/2+\ell/2}\le 2^{\ell/2}\le 2^{\frac{j+\ell}4}\le 2^{\frac{j+k+\ell}4}.
\end{align*}

\end{itemize}
\item If $(Q_A(\xi),Q_B(\xi))=(\xi_1\xi_2,\xi_1\xi_3+\xi_2^2),$ the computations are then as follows:
\begin{itemize}
    \item If $2k\ge j+\ell-100,$ then
$\frac j2=\frac j3+\frac j6\le\frac j3+\frac{2k-\ell+100}6\le 17+\frac{j+k}3\le 17+\frac{j+k+\ell}3$. Freezing the first variable, we obtain
$$2^{-j/2}|d\Phi_{1,\eta_1}(\xi,\eta_1,\eta_2)|^{-1/2}=
2^{-j/2}\big(|\xi_1-\eta_1|\big)^{-1}\sim2^{j/2}\lesssim 2^{\frac{j+k+\ell}3}.$$
    \item On the other hand, if $2k< j+\ell-100,$ then $2^{-2k}\gg2^{-j-\ell}$, hence, freezing the third variable now,
\begin{align*}
    2^{-\ell/2}|d\Phi_{3,\eta_3}(\xi,\eta_1,\eta_2)|^{-1/2}=&2^{-\ell/2}\big(2|\xi_2-\eta_2|^2-(\xi_1-\eta_1)(\xi_3-\eta_3)\big)^{-1/2}\\&\sim2^{-\ell/2}\big( 2\cdot 2^{-2k}-2^{-j-\ell}\big)^{-1/2}\\
    &\lesssim2^{k-\ell/2}\le 2^{k/3+2k/3-\ell/2}\\
    &\lesssim 2^{k/3+(j+\ell)/3-\ell/2}\lesssim2^{\frac{j+k+\ell}3}.
\end{align*}
\end{itemize}
\end{enumerate}
\hfill\qed

\bigskip

\thispagestyle{empty}

\renewcommand{\refname}{References}

\end{document}